\documentclass[a4paper,12pt]{amsart}
\usepackage{amssymb}

\DeclareMathSymbol{\subsetneqq}{\mathbin}{AMSb}{36}

\newcommand{\R}{\mathbb{R}}
\newcommand{\Z}{\mathbb{Z}}
\newcommand{\N}{\mathbb{N}}

\newcommand{\C}{\mathbb{C}}

\newcommand{\dint}{\displaystyle\int}

\newcommand{\dlim}{\displaystyle\lim}

\newcommand{\dliminf}{\displaystyle\liminf}

\numberwithin{equation}{section}
\newtheorem{thm}{Theorem}[section]
\newtheorem{prop}[thm]{Proposition}
\newtheorem{lemma}[thm]{Lemma}
\newtheorem{cor}[thm]{Corollary}
\newtheorem{defi}[thm]{Definition}

\theoremstyle{remark}
\newtheorem{rem}[thm]{Remark}

\author{J.~Colliander}
\address{Department of Mathematics, University of Toronto, Toronto, Ontario, M5S 2E4 Canada.}
\email {\it colliand@math.toronto.edu}
\thanks{J.C. is partially  supported  by N.S.E.R.C. Grant RGPIN 250233-07.}

\author{S.~Ibrahim}
\address{Department of Mathematics and Statistics, Arizona State University, P.O. Box 871804
Tempe, AZ   85287-1804, USA.}
\email {\it ibrahim@math.asu.edu}
\thanks{{S. I. is grateful to the Departments of Mathematics and Statistics at McMaster University and Arizona State University}}
\author{M. Majdoub}
\address{Facult\'e des Sciences de Tunis, D\'epartement
de Math\'ematiques\\ Campus universitaire 1060, Tunis, Tunisia.}
\email{\it mohamed.majdoub@fst.rnu.tn}
\thanks{M.M is partially supported by the {\sf Laboratory of PDE and applications}
of Faculty of Sciences, Tunis, Tunisia}
\author{N. Masmoudi}
\address{Department of Mathematics, Courant Institute of Mathematical
Sciences New York University, 251 Mercer\\ St. New York, NY 10012,
U.S.A.}
\email{\it masmoudi@cims.nyu.edu}
\thanks{N. M is partially supported by an NSF Grant DMS-0703145}

\title{Energy Critical NLS in two space dimensions}
\date{\today}
\begin{document}
\begin{abstract}
We investigate the initial value problem for a defocusing nonlinear Schr\"odinger
equation with exponential nonlinearity
$$
i\partial_t u+\Delta
u=u\big(e^{4\pi|u|^2}-1\big)\quad\mbox{in}\quad \R_t\times\R_x^2.
$$
We identify subcritical, critical and supercritical regimes in the
energy space. We establish global well-posedness in the subcritical
and critical regimes. Well-posedness fails to hold in the supercritical case.
\end{abstract}

%@@@@@@@@@@@@@@@@@@@@@@@@@@@@@@@@@@@@%@@@@@@@@@@@@@@@@@@@@@@@@@@@@@@@@@@@@%@@@@@@@@@@@@@@


%@@@@@@@@@@@@@@@@@@@@@@@@@@@@@@@@@@@@%@@@@@@@@@@@@@@@@@@@@@@@@@@@@@@@@@@@@%@@@@@@@@@@@@@@
\maketitle

\tableofcontents
\vspace{ -1\baselineskip}

\eject
%@@@@@@@@@@@@@@@@@@@@@@@@@@@@@@@@@@@@%@@@@@@@@@@@@@@@@@@@@@@@@@@@@@@@@@@@@%@@@@@@@@@@@@@@@
%@@@@@@@@@@@@@@@@@@@@@@@@@@@@@@@@@@@@%@@@@@@@@@@@@@@@@@@@@@@@@@@@@@@@@@@@@%@@@@@@@@@@@@@@@
\section{What is the energy critical NLS equation on $\R^2$?}
\setcounter{equation}{0} \setcounter{equation}{0}

We consider the initial value problem for a defocusing nonlinear Schr\"odinger  equation
with exponential nonlinearity
\begin{equation}
\label{eq1}
\left\{
\begin{matrix}
i\partial_t u+\Delta u= f(u), & u:(-T_*, T^*) \times \R^2 \longmapsto \C \\
u(0, \cdot) = u_0 (\cdot )\in H^1 (\R^2) \\
\end{matrix}
\right.
\end{equation}
where
\begin{equation}
\label{eq2}
f(u)=u\big(e^{4\pi|u|^2}-1\big).
\end{equation}
Solutions of (\ref{eq1}) formally satisfy the conservation of mass and Hamiltonian

\begin{eqnarray}
\label{eq3}
M(u(t,\cdot))&:=&\|u(t,\cdot)\|_{L^2}^2\\
\nonumber    &=&M(u(0,\cdot)),
\end{eqnarray}

\begin{eqnarray}
\label{eq4} H(u(t,\cdot))&:=&\Big\|\nabla u(t,\cdot)\Big\|_{L^2}^2+
\frac{1}{4\pi}\Big\|e^{4\pi |u(t,\cdot)|^2}-1-4\pi|u(t,\cdot)|^2\Big\|_{L^1(\R^2)}\\
\nonumber &=&H(u(0,\cdot)).
\end{eqnarray}

We show that for initial data $u_0$ satisfying $H(u_0)\leq 1$ the
initial value problem is global-in-time well-posed. Well-posedness fails to
hold for data satisfying $H(u_0 )> 1$.  We compare our theory for
\eqref{eq1} with work on the $\dot{H}^1$ critical NLS initial value problem
on $\R^d$ with $d \geq 3.$ Similar ill-posedness results were also obtained for the nonlinear Klein-Gordon equation with exponential nonlinearity in \cite{IMM5}.

\subsection{$NLS_p (\R^d)$ and critical regularity for local well-posedness}

We introduce a family of equations and identify \eqref{eq1} as a natural extreme
 limit of the family with monomial (or polynomial) nonlinearities when the space dimension is 2. The monomial defocusing
 semilinear initial value problem
\begin{equation}
\label{NLSp}
\left\{
\begin{matrix}
i \partial_t u + \Delta u =  |u|^{p-1} u, & u: (-T_*, T^* ) \times \R^d \longmapsto \C \\
u(0,x) = u_0 (x)
\end{matrix}
\right.
\end{equation}
has solutions which also satisfy conservation of mass and Hamiltonian, where
\begin{equation}
\label{NLSpEnergy}
H_p(u(t,\cdot)):=\|\nabla u(t,\cdot)\|_{L^2}^2+
\int_{\R^d} \frac{1}{p+1} |u|^{p+1} (t, x) dx.
\end{equation}
We will sometimes refer to the initial value problem \eqref{NLSp}
with the notation $NLS_p (\R^d)$.

If $u$ solves \eqref{NLSp} then, for $\lambda > 0,
u^\lambda: (-T_* \lambda^2 , T^* \lambda^2) \times \R^d$ defined by
\begin{equation}
\label{dilationsymmetry}
u^\lambda ( t , x ) : = \lambda^{2/(1-p)} u ( \lambda^{-2} t , \lambda^{-1} x)
\end{equation}
also solves \eqref{NLSp}. It turns out that Banach spaces whose
norms are invariant under the dilation $u \longmapsto u^\lambda$
are relevant in the theory of the initial value problem
\eqref{NLSp}. Let $s_c = \frac{d}{2} - \frac{2}{p-1}$. Note that
for all $\lambda >0$ the $L^2$-based homogeneous $\dot{H}^{s_c}$
Sobolev norm is invariant under the the mapping $f(x) \longmapsto
\lambda^{-2/(p-1)} f(\lambda^{-1} x)$. Similarly, note that the Lebesgue
$L^{p_c} (\R^d_x )$ norm is invariant under the dilation symmetry for
$p_c = \frac{d}{2} (p-1).$ Unless otherwise stated, we
will restrict{\footnote{Global well-posedness for the defocusing
energy supercritical $NLS_p (\R^d)$ with $s_c = \frac{d}{2} - \frac{2}{p-1} > 1$

is an outstanding open problem.}} this discussion to problems where
dimension $d$ and the degree $p$ are constrained to give $ 0 \leq
s_c \leq 1$. It is now known ({\cite{CazWeis}, \cite{GV}, \cite{CCT1}) that \eqref{NLSp}
with $H^s$ initial data is locally well-posed for $s > s_c$ with
existence interval depending only upon $\| u_0 \|_{H^s}$, locally
well-posed for $s=s_c$ with existence interval depending upon
$e^{it \Delta } u_0$, and is ill-posed for $s< s_c$.  Based on
this complete trichotomy, it is natural to refer to $H^{s_c}$ as
the {\it{critical regularity}} for \eqref{NLSp}.

\subsection{Global well-posedness for $NLS_p (\R^d)$}

For the {\it{energy subcritical}} case, when $s_c < 1$, an iteration
of the local-in-time well-posedness result using the {\it{a priori}}
upper bound on $\| u(t) \|_{H^1}$ implied by the conservation laws
establishes global well-posedness for \eqref{NLSp} in $H^1$. It is
expected that the local-in-time $H^{s_c}$ solutions of \eqref{NLSp}
extend to global-in-time solutions. For certain choices of $p, d$ in
the energy subcritical case, there are results
(\cite{B98}, \cite{B99}, \cite{CKSTT:MRL}, \cite{Tzirakis}, \cite{CKSTT:CPAM}) which establish that $H^s$
initial data $u_0$ evolve into global-in-time solutions $u$ of
\eqref{NLSp} for $s \in (\tilde{s}_{p,d} , 1)$ with $s_c <
{\tilde{s}}_{p,d} < 1$ and $\tilde{s}_{p,d} $ close to $1$ and away
from $s_c$. For all problems with $0 \leq s_c < 1$, global
well-posedness in the scaling invariant space $H^{s_c}$ is unknown
but conjectured to hold.

For the {\it{energy critical}} case, when $s_c =1$, an iteration of
the local-in-time well-posedness theory fails to prove global
well-posedness. Since the local-in-time existence interval depends
upon absolute continuity properties of the linear evolution $e^{it
\Delta} u_0$ (and not upon the controlled norm $\| u(t) \|_{H^1}$),
the local theory does not directly globalize based on the
conservation laws. Nevertheless, based on new ideas of Bourgain in
\cite{B99} (see also
\cite{BourgainBook})}
 (which treated the radial case in
dimension 3) and a new interaction Morawetz inequality
\cite{CKSTT:CPAM}
 the energy critical case of
\eqref{NLSp} is now completely resolved
\cite{CKSTT:ANMATH, Visan, TaoVisan} : Finite energy initial data $u_0$ evolve
into global-in-time solutions $u$ with finite spacetime size $\| u
\|_{L^{\frac{2(2+d)}{d-2}}_{t,x}} < \infty$ and scatter.

\subsection{Energy criticality in two space dimensions}

The initial value problem $NLS_p (\R^2)$ is energy subcritical for
all $p > 1$. To identify an "energy critical" nonlinear
Schr\"odinger initial value problem on $\R^2$, it is thus natural
to consider problems with exponential nonlinearities. In this
paper, we establish local and global well-posedness for
\eqref{eq1} provided that $H(u_0)\leq 1$. The case where
$H(u_0)= 1$ is more subtle than the case where $H(u_0)< 1.$ We also establish that well-posedness fails to hold
on the set of initial data where $H(u_0)> 1$. Thus, we
establish a complete trichotomy analogous to the energy critical
cases of $NLS_p (\R^d)$ in dimensions $d \geq 3.$ Based on these
results, we argue that \eqref{eq1} should be viewed as
{\underline{the}} energy critical NLS problem on $\R^2$. Using
a new interaction Morawetz estimate, proved independently by Colliander-Grillakis-Tzirakis and Planchon-Vega
\cite{CGT, PV}, the scattering was recently shown in \cite{IMMN} for subcritical solutions of \eqref{eq1} (with $f(u)=u\big(e^{4\pi|u|^2}-1-4\pi|u|^2\big)$). This problem remains open when $H(u_0)=1$ due to the lack of uniform global estimates of the nonlinear term and the infinite speed of propagation.

\begin{rem}
\label{energysupercriticalremark}
The critical threshold for local and global well-posedness of \eqref{eq1}
is expressed in terms of the size of $H(u_0)$. In contrast, the critical
threshold for energy critical \eqref{NLSp} is expressed in terms of $\| u_0 \|_{H^1}$.
Positive results for data satisfying $H(u_0)> 1$ and other conditions may give
insights towards proving global well-posedness results for energy supercritical problems.
\end{rem}

\subsection{Statements of results}

We begin by formally defining our notion of criticality and
well-posedness for \eqref{eq1}. We then give precise statements of
the main results we obtain and make brief comments about the rest
of the paper.

\begin{defi}
\label{d1} The Cauchy problem associated to (\ref{eq1}) and with
initial data $u_0\in H^1(\R^2)$ is
said to be {\it subcritical} if
$$
H(u_0)<1.
$$
It is {\it critical} if $H(u_0)=1$ and
{\it supercritical} if $H(u_0)>1$.
\end{defi}

\begin{defi}
\label{d2}
 We say that the Cauchy problem associated to
(\ref{eq1}) is {\it locally well-posed} in $H^1(\R^2)$ if there
exist $E>0$ and a time $T=T(E)>0$ such that for every $u_0\in
B_E:=\{\;u_0\in H^1(\R^2);\quad \|\nabla u_0\|_{L^2}<E\;\}$ there
exists a unique (distributional) solution $u:
[-T,T]\times\R^2\longrightarrow\C$ to (\ref{eq1}) which is in the
space ${\mathcal C}([-T,T]; H^1_x)$, and such that the solution
map $u_0\longmapsto u$ is uniformly continuous from $B_E$ to
${\mathcal C}([-T,T]; H^1_x)$.
\end{defi}

A {\it priori,} one can estimate the nonlinear part of the energy
\eqref{eq4} using the following Moser-Trudinger type
inequalities (see \cite{Tan},  \cite{M}, \cite{Tru}).

\begin{prop}[\sf Moser-Trudinger Inequality]\quad\\
\label{p1}
Let $\alpha\in [0,4\pi)$. A constant $c_\alpha$ exists
such that
\begin{equation}
\label{eq5} \|\exp(\alpha |u|^2)-1\|_{L^1(\R^2)}\leq c_\alpha
\|u\|_{L^2(\R^2)}^2
\end{equation}
for all $u$ in $H^1(\R^2)$ such that $\|\nabla
u\|_{L^2(\R^2)}\leq1$. Moreover, if $\alpha\geq 4\pi$, then
(\ref{eq5}) is false.
\end{prop}

\begin{rem}
\label{rem} We point out that $\alpha=4\pi$ becomes admissible in
(\ref{eq5}) if we require $\|u\|_{H^1(\R^2)}\leq1$ rather than
$\|\nabla u\|_{L^2(\R^2)}\leq1$. Precisely, we have
$$\sup_{\|u\|_{H^1}\leq
1}\;\;\|\exp(4\pi |u|^2)-1\|_{L^1(\R^2)}<+\infty$$ and this is
false for $\alpha>4\pi$. See \cite{Ru} for more details.
\end{rem}
To establish an energy estimate, one has to consider the
nonlinearity as a source term in \eqref{eq1}, so we need to estimate
it in the $L^1_t(H^1_x)$ norm. To do so, we use (\ref{eq5}) combined
with the so-called Strichartz estimate.
\begin{prop}[\sf Strichartz estimates]\quad\\
\label{p3} Let $v_0$ be a function in $H^1(\R^2)$ and $F\in
L^1(\R, H^1(\R^2))$. Denote by $v$ the solution of the
inhomogeneous linear Schr\"odinger equation
$$
i\partial_t v+\Delta v=F
$$
with initial data $v(0,x)=v_0(x)$. \\
Then, a constant $C$ exists such that for any $T>0$ and any
admissible couple of Strichartz exponents $(q,r)$ i.e
$0\leq\frac{2}{q}=1-\frac{2}{r}<1$, we have

\begin{equation}
\label{eq7} \|v\|_{L^q([0,T],{\mathcal B}^{1}_{r,2}(\R^2))}\leq
C\left[\|v_0\|_{H^1(\R^2)}+\| F\|_{L^1([0,T],H^1(\R^2))}\right].
\end{equation}
\end{prop}
In particular, note that $(q,r)=(4,4)$ is an admissible Strichartz
couple and

\begin{equation}
{\mathcal B}^{1}_{4,2}(\R^2)\longrightarrow {\mathcal
C}^{1/2}(\R^2).
\end{equation}
Recall that, for $1\leq
p,\; q\leq\infty$ and $s\in\mathbb{R}$, the (inhomogeneous) Besov
norm $\|.\|_{{\mathcal B}_{p,q}^s(\mathbb{R}^2)}$ is defined by
$$ \|u\|_{{\mathcal
B}_{p,q}^s(\mathbb{R}^2)}:=\Big(\sum_{j\geq-1}2^{jqs}\|{\Delta}_j
u\|_{L^p}^q\Big)^{\frac{1}{q}} $$ with the usual modification when
$q=\infty$. $\Big(\Delta_j\Big)$ is a (inhomogeneous) dyadic
partition of unity.
\begin{rem}\quad\\
$\bullet$ The homogeneous Besov norm is defined in the same manner
using a homogeneous dyadic partition of unity
$\Big(\dot{\Delta}_j\Big)_{j\in\Z}$.\\
$\bullet$  The connection between Besov spaces and the usual
Sobolev and H\"older spaces is given by the following relations
$${\mathcal B}_{2,2}^s(\mathbb{R}^2)=H^s(\mathbb{R}^2),\quad
{\mathcal B}_{\infty,\infty}^s(\mathbb{R}^2)={\mathcal
C}^s(\mathbb{R}^2).$$
\end{rem}
We recall without proof the following properties of Besov spaces
(see \cite{Tr1}, \cite{Tr2} and \cite{Tr3}).
\begin{thm}[{\rm\sf Embedding result}]\quad\\
The following injection holds
$$ {\mathcal B}_{p,q}^s(\mathbb{R}^2)\hookrightarrow{\mathcal
B}_{p_1,q_1}^{s_1}(\mathbb{R}^2) $$ where $$ \left\{
\begin{array}{cllll}
&s-\frac{2}{p}=s_1-\frac{2}{p_1},\\\\
 &1\leq p\leq
p_1\leq\infty,\quad 1\leq q\leq q_1\leq\infty,\quad s,s_1\in
\mathbb{R}.
\end{array}
\right. $$
\end{thm}

The following estimate is an $L^\infty$ logarithmic inequality
which enables us to establish the link between $ \|e^{4\pi
|u|^2}-1\|_{L^1_T(L^2(\R^2))} $ and dispersion properties of
solutions of the linear Schr\"odinger equation.

\begin{prop}[\sf Log Estimate]\quad\\
\label{p2}
Let $\beta\in]0,1[$. For any $\lambda>\frac{1}{2\pi\beta}$ and
any $0<\mu\leq1$, a constant $C_{\lambda}>0$ exists such
that,
for any function $u\in H^1(\R^2)\cap{\mathcal C}^\beta(\R^2)$, we have
\begin{equation}
\label{eq6}
\|u\|^2_{L^\infty}\leq
\lambda\|u\|_\mu^2 \log(C_{\lambda} +
\frac{8^\beta\mu^{-\beta}\|u\|_{{\mathcal C}^{\beta}}}{\|u\|_\mu}),
\end{equation}
where we set

\begin{equation}
\label{eq6'}
\|u\|_\mu^2:=\|\nabla u\|_{L^2}^2+\mu^2\|u\|_{L^2}^2.
\end{equation}
\end{prop}
Recall that ${\mathcal C}^{\beta}(\R^2)$ denotes the space of
$\beta$-H\"older continuous functions endowed with the norm $$
\|u\|_{{\mathcal
C}^{\beta}(\R^2)}:=\|u\|_{L^\infty(\R^2)}+\sup_{x\neq
y}\frac{|u(x)-u(y)|}{|x-y|^{\beta}}. $$
We refer to  \cite{IMM} for the proof of this proposition and more
details. We just point out that the condition  $ \lambda>
\frac{1}{2\pi\beta}$ in (\ref{eq6}) is optimal.

Our first
statement describes a local well-posedness result when the initial
data is in the open unit ball of the homogeneous Sobolev space
${\dot H}^1(\R^2)$. The sign of the nonlinearity is irrelevant here. Consider the following equation :
\begin{eqnarray}
\label{sigma} i\partial_t u+\Delta u= \sigma f(u).
\end{eqnarray}
We have the following short time existence Theorem.
\begin{thm}
\label{t1} Let $\sigma\in\{-1,+1\}$ and $u_0\in H^1(\R^2)$ such
that $\|\nabla u_0\|_{ L^2(\R^2)}<1$. Then, there exists a
time $T>0$ and a unique solution to the equation (\ref{sigma}) in
the space ${\mathcal
C}_T(H^1(\R^2))$ with initial data $u_0$ .  \\
Moreover, $u\in L^4_T({\mathcal C}^{1/2}(\R^2))$ and satisfies,
for all $0\leq t < T$,   $ M(u(t,\cdot)) = M(u_0)$ and
$H(u(t,\cdot))=H(u_0)$.
\end{thm}
The proof of this Theorem is similar to Theorem 1.8 in
\cite{IMM1}. It is based on the combination of the three {\it a priori} estimates
given by the above propositions. We derive the local
well-posedness using a classical fixed point argument.

\begin{rem}\label{r5}
In \cite{IMM6} a weak well-posedness result was proved without any
restriction on the size of the initial data. More precisely, it is
shown that the solution map is only continuous, while Theorem \ref{t1}
says that it is uniformly continuous when $\|\nabla u_0\|_{
  L^2(\R^2)}<1$. Well-posedness results with merely continuous
dependence upon the initial data have also been obtained for the KdV
equation \cite{KappelerTopalov} using the completely integrable
machinery and for the cubic NLS on the line \cite{CCT:NLS},
\cite{KochTataru} using PDE methods.
\end{rem}

\begin{rem}\label{r6}
In the defocusing case, the assumption $H(u_0)\leq 1$ in particular implies that
$\|\nabla u_0\|_{ L^2(\R^2)}<1$, and consequently we have the
short-time existence of solutions in both subcritical and critical
case. So it makes sense to investigate global existence in these
cases.
\end{rem}

As an immediate consequence of Theorem \ref{t1} we have the
following global existence result.
\begin{thm}[{\sf Subcritical case}]\quad\\
\label{tsub} Assume that $H(u_0)< 1$; then the defocusing problem
\eqref{eq1} has a unique global solution $u$ in the class
$$
{\mathcal C}(\R, H^1(\R^2)).
$$
Moreover, $u\in L^4_{loc}(\R,\;{\mathcal C}^{1/2}(\R^2))$ and
satisfies the conservation laws \eqref{eq3} and \eqref{eq4}.
\end{thm}
The reason behind Definition \ref{d1} is the following: If $u$ denotes
the solution given by Theorem \ref{t1}, where $T^*<\infty$ is the
largest time of existence, then the conservation of the total energy
gives us, in the subcritical setting, a uniform bound of $\|\nabla
u(t,\cdot)\|_{L^2(\R^2)}$  away from $1$, and therefore the solution
can be continued in time. In contrast, for the critical case, we
lose this uniform control and the total energy can be concentrated
in the $\|\nabla u(t,\cdot)\|_{L^2(\R^2)}$ part. By using a
localization result due to Nakanishi (see Lemma 6.2 in \cite{Nak}),
we show that such concentration cannot hold in the critical case and
therefore we have the following theorem:
\begin{thm}[{\sf Critical case}]\quad\\
\label{tc} Assume that $H(u_0)=1$; then the problem
\eqref{eq1} has a unique global solution $u$ in the class
$$
{\mathcal C}(\R, H^1(\R^2)).
$$
Moreover, $u\in L^4_{loc}(\R,\;{\mathcal C}^{1/2}(\R^2))$ and
satisfies the conservation laws \eqref{eq3} and \eqref{eq4}.
\end{thm}

\begin{rem}
Recently in \cite{IMMN} the scattering was established in the subcritical case using a new estimate obtained independently in \cite{CGT, PV}.
\end{rem}

When the initial data are more regular, we can easily prove that
the solution remains regular. More precisely, we have the
following theorem:
\begin{thm}
\label{treg} Assume that $u_0\in H^{s}(\R^2)$ with $s>1$ and
$\|\nabla u_0\|_{L^2(\R^2)}<1$. Then, the solution $u$ given in
Theorem \ref{t1} is in the space ${\mathcal C}_T(H^{s}(\R^2))$.
\end{thm}

\begin{rem}
In fact, the local well-posedness holds in $H^s$ for $s>1$ without
any assumption on the size of the initial data.
\end{rem}
The last result in this paper concerns the supercritical case.

\begin{thm}[{\sf The supercritical case}]\quad\\
\label{illposedness}
There exist sequences of initial data $u_k(0)$ and $v_k(0)$ bounded in $H^1$ and satisfying
$$
\dliminf_{k\to\infty}\,H(u_k(0))>1,\quad \dliminf_{k\to\infty}\,H(v_k(0))>1,
$$
with
$$
\lim_{k\rightarrow+\infty}\| u_k(0) - v_k(0) \|_{H^1}= 0,
$$
but there exists a sequence of times $t_k >0$ with $t_k \rightarrow0$ and
$$
\dliminf_{k\to\infty}\,\| \nabla (u_k (t_k) - v_k (t_k)) \|_{L^2} \gtrsim 1.
$$
\end{thm}

\begin{rem}
The sequences of initial data constructed in Theorem
\ref{illposedness} do not have bounded Hamiltonians. Indeed, their
potential parts are huge. Unlike for the Klein-Gordon where the speed
of propagation is finite see \cite{IMM5,IMM6}, we were unable to prove
the above result for slightly supercritical data.
\end{rem}

This class of two-dimensional problems with exponential growth
nonlinearities has been studied, for small Cauchy data, by Nakamura and Ozawa in
\cite{NO}. They proved global well-posedness and scattering.\\

%@@@@@@@@@@@@@@@@@@@@@@@@@@@@@@@@@@@@%@@@@@@@@@@@@@@@@@@@@@@@@@@@@@@@@@@@@%@@@@@@@@@@

\noindent{\sf Notation.} Let $T$ be a positive real number. We
denote by $\mbox{{\bf X}}(T)$ the Banach space defined by
$$
\mbox{{\bf X}}(T)={\mathcal C}_T(H^1(\R^2))\cap L^4_T({\mathcal
C}^{1/2}(\R^2)),
$$
and endowed with the norm
$$
\|u\|_{T}:=\sup_{t\in [0,T]}\big(\|u(t,\cdot)\|_{L^2}+ \|\nabla
u(t,\cdot)\|_{L^2}\big)+\|u\|_{L^4_T({\mathcal C}^{1/2})}.
$$
Here and below ${\mathcal C}_T( X )$ denotes ${\mathcal C }
([0,T); X )$ and $L^p_T(X)$ denotes $L^p([0,T) ; X)$.

If $A$ and $B$ are nonnegative quantities, we use $A\lesssim B$ to
denote $A\leq C B$  for some positive universal constant $C$, and
$A\thickapprox B$
to denote the estimate $A\lesssim B\lesssim A$.\\
For every positive real number $R$, $B(R)$ denotes the ball in
$\R^2$ centered at the origin and with radius $R$.

%@@@@@@@@@@@@@@@@@@@@@@@@@@@@@@@@@@@@%@@@@@@@@@@@@@@@@@@@@@@@@@@@@@@@@@@@@%@@@@@@@@@@@@@@@@

\section{Local well-posedness}
This section is devoted to the proof of Theorem \ref{t1} about
local existence. We begin with the following Lemma which summarizes some
of the properties of the exponential nonlinearity.

\begin{lemma}[Nonlinear Inhomogeneous Estimate]
\label{l1} Let $f$ be the function given by (\ref{eq2}), $T>0$ and
$0\leq A<1$. There exists $0<\gamma=\gamma(A)<3$ such that for any
two functions $U_1$ and $U_2$ in $\mbox{{\bf X}}(T)$ satisfying
the following

\begin{eqnarray}
\label{eq11}
\sup_{t\in[0,T]}\|\nabla U_j(t,\cdot)\|_{L^2}\leq A,
\end{eqnarray}
we have

\begin{eqnarray}
\label{eq12'}\quad\quad
\|f(U_1)-f(U_2)\|_{L^1_T(H^1(\R^2))}&\lesssim& \|U_1-U_2\|_{T}
\Big\{T^{\frac{3}{4}}\sum_{j=1,2}\|U_i\|_{T}^3
\\
 \nonumber
&+& T^{\frac{3-\gamma}{4}}\sum_{j=1,2}
\Big\|\frac{U_j}{A}\Big\|_{T}^\gamma\Big\}
\end{eqnarray}
\end{lemma}
%@@@@@@@@@@@@@@@@@@@@@@@@@@@@@@@@@@@@%@@@@@@@@@@@@@@@@@@@@@@@@@@@@@@@@@@@@%@@@@@@@@@@@@@@@@
\begin{proof}[Proof of Lemma \ref{l1}]
Let us identify $f$ with the ${\mathcal C}^\infty$ function defined on $\R^2$ and
denote by $D f$ the $\R^2$ derivative of the identified function. Then using the
mean value theorem and the convexity of the exponential function, we derive the
following  properties:
\begin{eqnarray}
\nonumber |f(z_1)-f(z_2)|&\lesssim&
|z_1-z_2|\sum_{j=1,2}\left(e^{4\pi|z_i|^2}-1 +|z_j|^2
e^{4\pi|z_j|^2}\right),
\end{eqnarray}
and

\begin{eqnarray}
\nonumber |(Df)(z_1)-(Df)(z_2)|\lesssim
|z_1-z_2|\sum_{j=1,2}\left(|z_i|e^{4\pi|z_j|^2} +
|z_j|^3e^{4\pi|z_j|^2}\right)
\end{eqnarray}

Therefore,  for any positive real number $\varepsilon$ there
exists a positive constant $C_{\varepsilon}$ such that
\begin{eqnarray}
\label{eq12} |f(z_1)-f(z_2)|\leq C_{\varepsilon} |z_1-z_2|
\Big\{e^{4\pi(1+\varepsilon)|z_1|^2}-1+
e^{4\pi(1+\varepsilon)|z_2|^2}-1\Big\},
\end{eqnarray}
and
\begin{eqnarray}
\label{eq13} |(Df)(z_1)-(Df)(z_2)|\leq C_{\varepsilon}  |z_1-z_2|
\sum_{i=1,2}\left(|z_i|+ e^{4\pi(1+\varepsilon)|z_i|^2}-1\right).
\end{eqnarray}
Now we estimate $\|f(U_1)-f(U_2)\|_{L^1_T(L^2(\R^2))}$. Applying the
H\"older inequality and using (\ref{eq12}) we infer
\begin{eqnarray}
\nonumber
\|f(U_1)(t,\cdot)-f(U_2)(t,\cdot)\|_{L^1_T(L^2(\R^2))}\leq
C_{\varepsilon}
\|U_1-U_2\|_{L^4_T(L^4)}\sum_{j=1,2}
\|e^{4\pi(1+\varepsilon)|U_1(t,\cdot)|^2}-1\|_{L^{4/3}_T(L^4)}.
\end{eqnarray}
Applying H\"older inequality again, we obtain
\begin{equation}
%\noindent\;\;\;\quad\quad 
\Big\|e^{4\pi(1+\varepsilon)|U_j(t,\cdot)|^2}-1\Big\|_{L^{4/3}_T (L^4_x)} \leq
\Big\| e^{3\pi(1+\varepsilon)
\|U_j(t,\cdot)\|_{L^\infty_x}^2}\Big\|_{L^{4/3}_T}
\|e^{4\pi(1+\varepsilon)|U_j(t,\cdot)|^2}-1\|_{L^\infty_T (L^1_x)}^{\frac{1}{4}}
\label{eq144}.
\end{equation}
Thanks to the Moser-Trudinger inequality (\ref{eq5}) and the Log estimate
(\ref{eq6}) we get

\begin{eqnarray}
\label{eq14} \quad\quad\quad
\|e^{4\pi(1+\varepsilon)|U_j(t,\cdot)|^2}-1\|_{L^1}
\leq
C_{4\pi(1+\varepsilon)A^2}
\|U_j(t,\cdot)\|_{L^2}^{2},
\end{eqnarray}
\begin{eqnarray}
\label{eq15}
 e^{3\pi(1+\varepsilon)
\|U_j(t,\cdot)\|_{L^\infty_x}^2}\lesssim \left(e^3+\frac{\|U_j(t,\cdot)\|_{{\mathcal
C}^{1/2}}}{{A^{\prime}}}\right)^{\gamma},
\end{eqnarray}
where we set
$$
{A^{\prime}}^2:=A^2+\max_i\sup_{t\in
[0,T]}\mu^2\|U_j(t,\cdot)\|_{L^2}^2\quad\mbox{and}\quad \gamma:
=3\pi\lambda(1+\varepsilon){A^{\prime}}^2,
$$
and $0<\mu\leq 1$ is chosen such that $A^{\prime}<1$. Remember
that $ C_{4\pi(1+\varepsilon)A^2}$ is given by Proposition
\ref{p1}. It is important to note that estimate (\ref{eq14}) is
true as long as the parameter $\varepsilon$ is such that
$(1+\varepsilon)A^2<1$. Now, inserting this back into
\eqref{eq144}, and integrating with respect to time, we obtain

\begin{eqnarray}
\nonumber \|(f(U_1)-f(U_2))(t,\cdot)\|_{L^1_T(L^2(\R^2))}\lesssim
\|U_1-U_2\|_{L^4_T(L^4)}
\sum_{j=1,2}\Big\|e^3+\frac{\|U_i\|_{{\mathcal
C}^{1/2}}}{{A^{\prime}}}\Big\|_{L^{\frac{4}{3}\gamma}_T}^{\gamma}\|U_j\|_{L^\infty_T(L^2)}^{1/2}.
\end{eqnarray}

Now we estimate
$\|f(U_1)-f(U_2)\|_{L^1([0,T], {\dot H}^1(\R^2))}$. We write

\begin{eqnarray*}
D(f(U_1)-f(U_2))&=& [(Df)(U_1)-(Df)(U_2)]DU_1+ Df(U_2)D(U_1-U_2)\\
&:=&(I)+(II).
\end{eqnarray*}
To estimate $(I)$ we use (\ref{eq13}). Hence for any $\varepsilon>0$
we have

\begin{eqnarray}
\nonumber |(Df)(U_1)-(Df)(U_2)|\leq C_{\varepsilon} |U_1-U_2|
\sum_{j=1,2}\left(e^{4\pi(1+\varepsilon)|U_j|^2}-1+|U_j|\right) ,
\end{eqnarray}
and therefore

\begin{eqnarray*}
|(I)|&\lesssim& |U_1-U_2| |DU_1|\sum_{j=1,2}|U_j|\\
     &+&|U_1-U_2||DU_1|\sum_{j=1,2}\left(e^{4\pi(1+\varepsilon)|U_j|^2}-1\right).
\end{eqnarray*}
Applying H\"older inequality we infer

\begin{eqnarray*}
\|(I)\|_{L^2(\R^2)}&\lesssim&
\|U_1-U_2\|_{L^8(\R^2)}\|DU_1\|_{L^4(\R^2)}\sum_{j=1,2}
\|U_j\|_{L^8(\R^2)}\\
&+&\|U_1-U_2\|_{L^{4(1+\frac{1}{\varepsilon})}(\R^2))}
\|DU_1\|_{L^4(\R^2)}\sum_{j=1,2}
\|e^{4\pi(1+\varepsilon)|U_j|^2}-1\|_{L^{4(1+\varepsilon)}}.
\end{eqnarray*}
Using (\ref{eq14}) and integrating with respect to time we deduce that
\begin{eqnarray*}
\|(I)\|_{L^1_T(L^2)}&\lesssim&
T^{\frac{3}{4}}\|U_1-U_2\|_{L^\infty_T(L^8)}
\|DU_1\|_{L^4([0,T]\times\R^2)}
\sum_{j=1,2}\|U_j\|_{L^\infty_T(L^8)}\\
&+&\|U_1-U_2\|_{L^\infty_T(L^{4(1+\frac{1}{\varepsilon})})}
\|DU_1\|_{L^4([0,T]\times\R^2)} \sum_{j=1,2}\Big\|e^3+
\frac{\|U_j(t,\cdot)\|_{{\mathcal
C}^{\frac{1}{2}}}}{{A^{\prime}}}\Big\|_{L^{\frac{4}{3}\gamma}_T}^{\gamma}.
\end{eqnarray*}

 To estimate the
term $(II)$, we use (\ref{eq5}) with $U_1=0$. So thanks to the
H\"older inequality we get

\begin{eqnarray*}
 \|(II)\|_{L^2(\R^2)}&\lesssim&
\|D(U_1-U_2)\|_{L^4(\R^2)}\|U_2\|_{L^8(\R^2)}^2\\
&+&\|D(U_1-U_2)\|_{L^4(\R^2)}\|U_2\|_{L^{4(1+\frac{1}{\varepsilon})}(\R^2)}
\sum_{j=1,2}\|e^{4\pi(1+\varepsilon)|U_j|^2}-1\|_{L^{4(1+\varepsilon)}}.
\end{eqnarray*}
Then we proceed exactly as we did for term $(I)$.\\
Now since $A<1$, we can choose the parameter $\mu$ such that
${A^{\prime}}<1$. Then we chose $\varepsilon>0$ small enough and
$\lambda>\frac{1}{\pi}$ and close to $\frac{1}{\pi}$ such that
$\gamma<3$. Applying H\"older inequality (with respect to time) in
the above inequality and in (\ref{eq14}), we deduce (\ref{eq12'})
as desired.
\end{proof}
%@@@@@@@@@@@@@@@@@@@@@@@@@@@@@@@@@@@@%@@@@@@@@@@@@@@@@@@@@@@@@@@@@@@@@@@@@%@@@@@@@@@@@@@@@@
\begin{proof}[Proof of Theorem \ref{t1}] The proof is divided into
two steps.\\

\noindent{\bf First step:} {\sf Local existence.}\\ Let $v_0$ be the
solution of the free Schr\"odinger equation with $u_0$ as the Cauchy
data. Namely,

\begin{eqnarray}
\label{eq18}
i\partial_t v_0+\Delta v_0=0\\
\nonumber
v_0(0,x)=u_0.
\end{eqnarray}

For any positive real numbers $T$ and $\delta$, denote by
${\mathcal E}_T(\delta)$ the closed ball in $\mbox{{\bf X}}(T)$ of
radius $\delta$ and centered at the origin. On the ball ${\mathcal
E}_T(\delta) $, define the map $\Phi$ by\\

\begin{eqnarray}
\label{eq16'} v\longmapsto\Phi(v):=\tilde{v},
\end{eqnarray}
where

\begin{eqnarray}
\label{eq17}
i\partial_t \tilde{v}+\Delta \tilde{v}=(v+v_0)\big(e^{4\pi|v+v_0|^2}-1\big),\quad
\tilde{v}(0,x)=0.
\end{eqnarray}

Now the problem is to show that, if $\delta$ and $T$ are small
enough, the map $\Phi$ is well defined from ${\mathcal
E}_T(\delta)$ into itself and it is a contraction.

In order to
show that the map $\Phi$ is well defined, we need to estimate the
term
$\|(v+v_0)\big(\;e^{4\pi|v+v_0|^2}-1\big)\|_{L^1_T(H^1)}$.\\
Let $U_1:=v+v_0$. Obviously, $U_1\in \mbox{{\bf X}}(T)$. Moreover,
since
$$
\|\nabla v_0(t,\cdot)\|_{L^2}=\|\nabla u_0\|_{L^2}
$$
is conserved along time, and $\|\nabla u_0\|_{L^2}<1$, then
the hypothesis (\ref{eq11}) of Lemma \ref{l1} is satisfied. Now
taking $U_2=0$, applying (\ref{eq12'}) and choosing $\delta$ and $T$
small enough show that $\Phi$ is well defined.
We do similarly for the contraction.
\end{proof}

\noindent{\bf Second step:} {\sf Uniqueness in the energy space.}\\
The uniqueness in the energy space is a straightforward consequence
of the following lemma and Theorem \ref{t1}. Note that uniqueness in
 ${\bf X} (T)$ follows from the contraction argument. Here we are noting
 the stronger statement that uniqueness holds in a larger space.

\begin{lemma}
Let $\delta$ be a positive real number and $u_0\in H^1(\R^2)$ such
that $\|\nabla u_0\|_{L^2}<1$. If $u\in {\mathcal C}([0,T],
H^1(\R^2))$ is a solution of (\ref{eq1})-(\ref{eq2}) on $[0,T]$,
then there exists a time $0<T_\delta\leq T$ such that
$u\;;\;\nabla u\in L^4([0,T_\delta], L^4(\R^2))$ and

$$
\| u\|_{L^4([0,T_\delta]\times \R^2)}+ \|\nabla
u\|_{L^4([0,T_\delta]\times \R^2)}\leq\delta.
$$
\end{lemma}

\begin{proof}
Fix $a>1$ such that
\begin{equation}
\label{a:choice} a\|\nabla u_0\|_{L^2}^2<1.
\end{equation}
Then choose $\varepsilon>0$ such that

\begin{eqnarray}
\label{epsilon:a} (1+\varepsilon)^2\;a\;\|\nabla
u_0\|_{L^2}^2<1\quad\mbox{and}\quad
4\frac{(1+\varepsilon)^2}{\varepsilon}\frac{a}{a-1}\varepsilon^2\leq
1.
\end{eqnarray}
Denote by $V:=u-v_0$ with $v_0:=e^{it\Delta}u_0$. Note that $V$
satisfies

\begin{eqnarray}
\nonumber
i\partial_t V+\Delta V= -(V+v_0)\big(e^{4\pi|V+v_0|^2}-1\big).
\end{eqnarray}
According to the Strichartz inequalities, to prove that $V$ and
$\nabla V$ are in $L^4_{t,x}$ it is sufficient to estimate
$\nabla^j\Big[(V+v_0)\big(e^{4\pi|V+v_0|^2}-1\big)\Big]$ in
the dual Strichartz norm $\|\cdot\|_{L^{\frac{4}{3}}}$ with $j=0,1$.\\
By continuity of $t\rightarrow V(t,\cdot)$, one can choose a
time $0<T_1\leq T$ such that

$$
\sup_{[0,T_1]}\|V(t,\cdot)\|_{H^1}\leq\varepsilon.
$$
Moreover, observe that

$$
|V+v_0|^2\leq a |v_0|^2 + \frac{a}{a-1} |V|^2,
$$

$$
e^{x+y}-1=(e^x -1)(e^y -1)+(e^x -1)+(e^y -1),
$$
and
$$
xe^{x}\leq \frac{e^{(1+\varepsilon)x}-1}{\varepsilon}.
$$
We will only estimate the term with derivative, the other case is easier.

\begin{eqnarray*}
\nonumber |\nabla\big[(V+v_0)\big(e^{4\pi|V+v_0|^2}-1\big)\big]|
&\leq& |\nabla(V+v_0)\big(e^{4\pi|V+v_0|^2}-1\big)|\\&+&
|\nabla(V+v_0)|V+v_0|^2e^{4\pi|V+v_0|^2}|\\&\leq&
|\nabla(V+v_0)\big(e^{4\pi|V+v_0|^2}-1\big)|\\&+&
|\nabla(V+v_0)\big(e^{4\pi(1+\varepsilon)|V+v_0|^2}-1\big)|.
\nonumber
\end{eqnarray*}
Hence we only need to estimate
$\|\nabla(V+v_0)
\big(e^{4\pi(1+\varepsilon)|V+v_0|^2}-1\big)\|_{L^{\frac{4}{3}}_{t,x}}$.
 Applying the H\"older inequality we obtain

\begin{eqnarray}
\nonumber
\|\nabla(V+v_0)
\big(e^{4\pi(1+\varepsilon)|V+v_0|^2}-1\big)\|_{L^{\frac{4}{3}}_{t,x}}
\leq \|\nabla(V+v_0)\|_{L^\infty_tL^2_x}
\|e^{4\pi(1+\varepsilon)|V+v_0|^2}-1\|_{L^{\frac{4}{3}}_tL^{4}_x}.
\end{eqnarray}
Using the above observations we need to estimate the following three terms

$$
{\mathcal I}_1(t):=\|\big(e^{4\pi(1+\varepsilon)a|v_0|^2}-1\big)
\big(e^{4\pi(1+\varepsilon)\frac{a}{a-1}|V|^2}-1\big)\|_{L^{\frac{4}{3}}_tL^{4}_x}
$$

$$
{\mathcal I}_2(t):=\|e^{4\pi(1+\varepsilon)a|v_0|^2}-1\|_{L^{\frac{4}{3}}_tL^{4}_x},
$$
and

$$
{\mathcal I}_3(t):=\|e^{4\pi(1+\varepsilon)
\frac{a}{a-1}|V|^2}-1\|_{L^{\frac{4}{3}}_tL^{4}_x}.
$$
Applying  H\"older inequality we obtain

$$
{\mathcal I}_1(t)\leq
\|e^{4\pi(1+\varepsilon)a|v_0|^2}-1\|_{L^{\frac{4}{3}}_tL^{4(1+\varepsilon)}_x}
\|e^{4\pi(1+\varepsilon)
\frac{a}{a-1}|V|^2}-1\|_{L^\infty_tL^{4\frac{1+\varepsilon}{\varepsilon}}_x}.
$$
Now the choice of the parameters $\varepsilon$ and $a$ satisfying
 (\ref{a:choice})-(\ref{epsilon:a}) insures that

$$
\|e^{4\pi(1+\varepsilon)
\frac{a}{a-1}|V|^2}-1\|_{L^\infty([0,T_1],L^{4\frac{1+\varepsilon}
{\varepsilon}})}\leq C(\varepsilon,a).
$$
Also

\begin{eqnarray*}
\|e^{4\pi(1+\varepsilon)a|v_0|^2}-1\|_{L^{\frac{4}{3}}([0,T_1],L^{4(1+\varepsilon)})}
&\leq& C(\varepsilon,a)\|e^{\pi(1+\varepsilon)a|v_0(t,\cdot)|_{L^\infty}^2}\|_{L^4[0,T_1]}^3\\
&\leq&C(\varepsilon,a,T_1).
\end{eqnarray*}
Note that $\dlim_{S\rightarrow0}C(\varepsilon,a,S)=0$, hence
choosing $T_1$ small enough we derive the desired estimate. The
other terms can be estimated in a similar way. We omit the details
here.
\end{proof}
%@@@@@@@@@@@@@@@@@@@@@@@@@@@@@@@@@@@@%@@@@@@@@@@@@@@@@@@@@@@@@@@@@@@@@@@@@%@@@@@@@@

\section{Global well-posedness}
In this section, we start with a remark about the time of local
existence. Then we show that the solutions emerging from the
subcritical regime in the energy space extend globally in time by a
rather simple argument. The more difficult critical case is then
treated with a nonconcentration argument.
\begin{rem}
\label{rem-exi} In Theorem  \ref{t1}, the time of existence $T$
depends on $u_0$. However, in the case $\|\nabla u_0\|_{
L^2(\R^2)}^2 < 1-\eta$, this time of
existence depends only on $\eta$ and $\| u_0\|_{ L^2(\R^2)}^2 $.
\end{rem}
%@@@@@@@@@@@@@@@@@@@@@@@@@@@@@@@@@@@@%@@@@@@@@@@@@@@@@@@@@@@@@@@@@@@@@@@@@%@@@@@@@@

\subsection{Subcritical Case}
\noindent Recall that in the subcritical setting we have
$H(u_0)<1$. Since the assumption $H(u_0)<1$
particularly implies that
$$
\|\nabla u_0\|_{L^2}<1,
$$
it follows that the equation \eqref{eq1} has a unique maximal
solution $u$ in the space $\mbox{{\bf X}}(T^*)$ where $0<T^*\leq
+\infty$ is the lifespan of $u$. We want to show that $T^*=+\infty$
which means that our solution is global in time.

\begin{proof}[Proof of Theorem \ref{tsub}]
Assume that $T^*<+\infty$, then  by the conservation of the
Hamiltonian (identity (\ref{eq4})), we deduce that

$$
\sup_{t\in[0,T^*)}\|\nabla u(t,\cdot)\|_{L^2(\R^2)}\leq
H(u_0)<1.
$$
Now, let $0<s<T^*$ and consider the
following Cauchy problem
$$\left\{
\begin{array}{cllll}
i\partial_t v+\Delta v&=&f(v)\\
 v(s,x)&=&u(s,x)\in H^1(\R^2).
\end{array}
\right.
$$
A fixed point argument (as that used in the proof of Theorem
\ref{t1}) shows that there exists a nonnegative $\tau$ and an
unique solution $v$ to our problem on the interval $[s,\;s+\tau]$.
Notice that $\tau$ does not depend on $s$ (see Remark
\ref{rem-exi} above). Choosing $s$ close to $T^*$ such that $T^*
-s<\tau$ the solution $u$ can be continued beyond the time  $T^*$
which is a contradiction.\end{proof}
%@@@@@@@@@@@@@@@@@@@@@@@@@@@@@@@@@@@@%@@@@@@@@@@@@@@@@@@@@@@@@@@@@@@@@@@@@%@@@@@@@@

\subsection{Critical Case}
\noindent Now, we consider the case when $H(u_0)=1$, and we
want to prove a global existence result as in the subcritical
setting. \\

The situation here is more delicate than that in the subcritical
setting; in fact the arguments used there do not apply here. Let
us briefly explain the major difficulty.  Since
$H(u_0)=1$ and by the conservation identities (\ref{eq3})
and (\ref{eq4}), it is possible (at least formally) that a
concentration phenomena occurs, namely
$$
\limsup_{t\rightarrow T^*}\;\;\|\nabla u(t,\cdot)\|_{L^2}=1
$$
 where $u$ is the maximal solution and
$T^*<+\infty$ is the lifespan of $u$.
In such a case, we can not apply the previous argument to continue
the solution. The actual proof is based on proving that the concentration
phenomenon does not happen.

Arguing by contradiction we claim the following.
\begin{prop}
\label{p4}
Let $u$ be the maximal solution of \eqref{eq1} defined on
$[0,\;T^*)$, and assume that $T^*$ is finite. Then

\begin{equation}
\label{ec1}
\limsup_{t\rightarrow T^*}\;\;\|\nabla
u(t)\|_{L^2(\R^2)}=1,
\end{equation}
and
\begin{equation}
\label{ec2}
\limsup_{t\rightarrow T^*}\;\;\|u(t)\|_{L^4(\R^2)}=0.
\end{equation}
\end{prop}

\begin{proof}
\noindent Note that for all $0\leq t<T^*$ we have

\begin{eqnarray}
\|\nabla u(t)\|^2_{L^2(\R^2)}\leq H(u(t,\cdot)).
\nonumber
\end{eqnarray}
On the other hand, since the Hamiltonian is conserved, we have
\begin{eqnarray}
\limsup_{t\rightarrow T^*}\;\;\|\nabla u(t)\|_{L^2(\R^2)}\leq
1.
\nonumber
\end{eqnarray}
Assume that

\begin{eqnarray}
\limsup_{t\rightarrow T^*}\;\;\|\nabla
u(t)\|_{L^2(\R^2)}=L<1.
\nonumber
\end{eqnarray}
Then, a time $t_0$ exists such that
$0<t_0<T^*$ and

\begin{eqnarray}
t_0<t<T^*\Longrightarrow\|\nabla
u(t)\|_{L^2(\R^2)}\leq\frac{L+1}{2}.
\nonumber
\end{eqnarray}

Take a time $s$ such that $t_0<s<T^*<s+\tau$ where $\tau$ depends
only on $\frac{1-L}{2}.$ Using the local existence theory, we
can extend the solution $u$ after the time $T^*$ which is a
contradiction. This concludes the proof of (\ref{ec1}). \\To
establish (\ref{ec2}), it is sufficient to note that
\begin{eqnarray}
\nonumber
2\pi|u(t,x)|^4\leq\frac{e^{4\pi |u(t,x)^2}|-1}{4\pi}-|u(t,x)|^2
\end{eqnarray}
and then consider the Hamiltonian with \eqref{ec1}.
\end{proof}
To localize the concentration and get a contradiction, the proof
in the case of the nonlinear Klein-Gordon equation was crucially based on the property of finite
speed of propagation satisfied by the solutions (see \cite{IMM1}).
Here that property breaks down. Instead, we use the following
localization result due to Nakanishi (see Lemma 6.2 in
\cite{Nak}).
\begin{lemma}
\label{N} Let $u$ be a solution of \eqref{eq1} on $[0,T)$ with
$0<T\leq+\infty$ and suppose that $E:=H(u_0)+M(u_0)<\infty$. A
constant $C(E)$ exists such that, for any two  positive real numbers $R$ and $R'$ and for any $0<t<T$, the following holds:
\begin{eqnarray}
\label{massdistribution} \int_{B(R+R')}|u(t,x)|^2dx\geq
\int_{B(R)}|u_0(x)|^2dx-C(E)\frac{t}{R'}.
\end{eqnarray}
\end{lemma}
For the sake of completeness, we shall give the proof here.
\begin{proof}[Proof of Lemma 2.6 \cite{Nak}]
Let $d_R(x):=d(x,B(R))$ be the distance from $x$ to the ball $B(R)$.
Obviously we have $|\nabla d_R(x)|\leq 1$. Define the cut-off
function
$$
\xi(x):=h\left(1-\frac{d_R(x)}{R'}\right)
$$
where $h$ is a smooth function such that $h(\tau)=1$ if $\tau\geq 1$
 and $h(\tau)=0$ if $\tau\leq 0$. Note that $\xi$ satisfies
$$
\xi(x)=1 \quad\hbox{  if  }\quad x\in B(R),\quad\xi(x)=0\quad \hbox{
if }\quad |x|\geq R+R'\quad\mbox{and}\quad
\|\nabla\xi(x)\|_{L^\infty}\lesssim 1/R'.$$
 Multiplying equation (\ref{eq1}) by $\xi^2 \bar{u}$,
 integrating on $\R^2$ and taking the imaginary part, we get the following identity
$$
\partial_t \|\xi u\|_{L^2}^2=4\;\mbox{Im}\left(\int_{\R^2}\xi\nabla\xi
u\nabla{\bar u}\;dx\right)\geq -\frac{C(E)}{R'}.
$$
This completes the proof of the Lemma.
\end{proof}

\begin{proof}[Proof of Theorem \ref{tc}] The proof of Theorem \ref{tc}
is now straightforward. Assuming that $T^*<+\infty$ and applying
H\"older inequality to the left hand side of
(\ref{massdistribution}), we infer

\begin{eqnarray*}
 \dint_{B(R)}|u_0(x)|^2dx-C(E)\frac{t}{R'}\leq
 \sqrt{\pi}\; (R+R')\;\|u(t)\|_{L^4(\R^2)}^2.
\end{eqnarray*}
Taking first the {\it limsup} as $t$ goes to $T^*$ and then $R'$ to
infinity we deduce that $u_0$ should be zero which leads to a
contradiction and therefore the proof is achieved.
\end{proof}

%@@@@@@@@@@@@@@@@@@@@@@@@@@@@@@@@@@@@%@@@@@@@@@@@@@@@@@@@@@@@@@@@@@@@@@@@@%

\section{Instability of supercritical solutions of NLS}
The aim of this section is to show that the Cauchy problem
\eqref{eq1} is ill-posed for certain data satisfying $H(u_0)>1$. A
typical example of supercritical data is the function $f_k$
defined by:
\begin{eqnarray*}
 f_k(x)&=&\; \left\{
\begin{array}{cllll}0 \quad&\mbox{if}&\quad
|x|\geq 1,\\\\
-\dfrac{\log|x|}{\sqrt{k\pi}} \quad&\mbox{if}&\quad e^{-k/2}\leq
|x|\leq 1
,\\\\
\sqrt{\frac{k}{4\pi}}\quad&\mbox{if}&\quad |x|\leq e^{-k/2}.
\end{array}
\right.
\end{eqnarray*}
These functions were introduced in \cite{M} to show the optimality
of the exponent $4\pi$ in Trudinger-Moser inequality (see also
\cite{Tan}).\\ An easy computation shows that $\|\nabla
f_k\|_{L^2(\R^2)}=1$. Since the sequence of functions $f_k$ is not smooth enough, we begin by
regularizing it in a way that preserves its ``shape" i.e. : Let $\chi$ be a smooth function
 such that $0\leq \chi\leq 1$ and
\begin{eqnarray*}
\chi(\tau)&=&\; \left\{
\begin{array}{cllll}0 \quad&\mbox{if}&\quad
\tau\leq 3/2,\\\\
1 \quad&\mbox{if}&\quad \tau\geq 2.
\end{array}
\right.
\end{eqnarray*}
For every integer $k\geq 1$, let
$\eta_k(x)=\chi(e^{k/2}|x|)\chi(e^{k/2}(1-|x|))$ and
$\tilde{f}_k=\eta_k\;f_k$. An easy computation show that, for all
$j\geq 0$, we have
$$
\|\eta_k^{(j)}\|_{L^\infty}\leq e^{jk/2}.
$$
For any nonnegative  $\alpha$ and $A>0$, denote by
$$
g_{\alpha, A,
k}(y):=\left(1+\frac{\alpha}{k}\right)\tilde{f}_k(y)\varphi\left(\frac{y}{\nu_k(A)}\right),
$$
where $\varphi$ is a cut-off function such that
$$
\mbox{ supp}(\varphi)\subset B(2),\quad \varphi=1
\quad\mbox{on}\quad B(1),\quad 0\leq \varphi\leq 1,
$$
and the following choice of the scale $\nu$
$$
\nu_k(A)=\exp\left(-\frac{\sqrt{k}}{A}\right)
$$
The cut-off function $\varphi$ is made to insure that the rescaled $g_{\alpha, A,
k}(\nu_k(A) x)$ has a finite $L^2$ norm.
Now, let $u$ solve the Cauchy problem
\begin{equation}
\label{NLSx} \left\{
\begin{matrix}
i \partial_t u + \Delta_x u =  f(u)\\\\
u(0,x) = g_{\alpha, A, k}(\nu_k(A) x).
\end{matrix}
\right.
\end{equation}
Define $v(t,\nu_k(A) x))=u(t,x)$. Then $v$ satisfies

\begin{equation}
\label{NLSy} \left\{
\begin{matrix}
i \partial_t v + \nu_k(A)^2\Delta_y v =  f(v)\\\\
v(0,y) = g_{\alpha, A, k}(y)
\end{matrix}
\right.
\end{equation}
For the sake of clarity, we omit the dependence of $u$ and $v$ upon the parameters $\alpha$, $k$ and $A$.
We begin by showing that the initial data is supercritical. 
\begin{lemma}
\label{sup-crit} There exists a positive constant $C_1$ such that for every $A>0$, we have

$$
\liminf_{k\rightarrow\infty}\;H(g_{\alpha, A,
k}(\nu_k(A)\cdot))\;\geq\;1+\frac{C_1}{\pi A^2}.
$$
\end{lemma}
\begin{proof}[Proof of Lemma \ref{sup-crit}]
For simplicity, we shall denote $g_{\alpha, A, k}$ by $g$ and
$\nu_k(A)$ by $\nu$. Recall that, by definition, we have
$$
g(y)=0\quad\mbox{if}\quad |y|\geq 2\nu\quad\mbox{or}\quad |y|\leq
\frac{3}{2}e^{- k/2}
$$
and
$$
g(y)=\left(1+\frac{\alpha}{k}\right)f_k(y)\quad\mbox{if}\quad
|y|\leq \nu\quad\mbox{and}\quad 2 e^{- k/2}\leq|y|\leq 1-2e^{- k/2}
$$

Remark that
$$
H(g)\geq\|\nabla g\|_{L^2}^2\geq (I)+ (II),
$$
where
\begin{eqnarray*}
(I)&=&\|\nabla g\|_{L^2(2 e^{- k/2}\leq|y|\leq
\nu)}^2\\\\
(II)&=&\|\nabla g\|_{L^2(\nu\leq|y|\leq 2\nu)}^2.
\end{eqnarray*}
On the set $\{2 e^{- k/2}\leq|y|\leq \nu\}$ we
have $g(y)=-(1+\alpha/k)\frac{\log|y|}{\sqrt{k\pi}}$ and thus
\begin{eqnarray}
\label{I}
(I)=1-\frac{2}{A\sqrt{k}}+\frac{2(\alpha-\log2)}{k}-\frac{4\alpha}{A
k^{3/2}}-\frac{4\alpha\log2}{k^2}.
\end{eqnarray}
For the second term, we write
\begin{eqnarray}
\label{II}
(II)= (a) + (b) + (c)
\end{eqnarray}
where
\begin{eqnarray*}
(a)&=&(1+\frac{\alpha}{k})^2\dint|\nabla f_k(y)|^2|\varphi(\frac{y}{\nu})|^2dy\\\\
(b)&=&(1+\frac{\alpha}{k})^2\nu^{-2}\dint|f_k(y)|^2|\nabla\varphi(y)|^2dy\\\\
(c)&=&2(1+\frac{\alpha}{k})^2\nu^{-1}\dint
f_k(y)\varphi(\frac{y}{\nu})\nabla
f_k(y)\cdot\nabla\varphi(\frac{y}{\nu})dy
\end{eqnarray*}
Clearly
\begin{eqnarray}
\label{a}
\indent(a)=
\frac{2}{k}(1+\frac{\alpha}{k})^2\left(\int_1^2\frac{|\varphi(r)|^2}{r}dr\right),
\end{eqnarray}
and
\begin{eqnarray*}
(b)=\frac{1}{\pi
k\nu^2}(1+\frac{\alpha}{k})^2\dint_{\nu\leq|y|\leq
2\nu}\log^2|y||\nabla\varphi(\frac{y}{\nu})|^2dy.
\end{eqnarray*}
But since, for
$k$ large, $\left(\log 2-\frac{\sqrt{k}}{A}\right)^2\leq \log^2|y|\leq \frac{k}{A^2}$,
we deduce that
\begin{eqnarray*}
\frac{1}{\pi k}(1+\frac{\alpha}{k})^2\left(\log
2-\frac{\sqrt{k}}{A}\right)^2\left(\dint_{1\leq|z|\leq
2}\;|\nabla\varphi(z)|^2dz\right)\leq(b)\leq\frac{C}{\pi A^2}(1+\alpha/k)^2,
\end{eqnarray*}
and therefore,

\begin{eqnarray}
 \label{b}
(1+\frac{\alpha}{k})^2\left(\frac{C_1}{\pi A^2}-\frac{2\log2
C_1}{\pi A\sqrt{k}}+\frac{C_1\log^2 2}{\pi k}\right)\leq(b)\leq\frac{C}{\pi A^2}(1+\alpha/k)^2.
\end{eqnarray}
The constant $C_1=\|\nabla\varphi\|_{L^2}^2$. For the last term, we simply
write
\begin{eqnarray}
\nonumber
(c)&=&\frac{2}{\pi k}(1+\frac{\alpha}{k})^2\dint_{1\leq|z|\leq
2}\;\left(\log\nu+\log|z|\right)\varphi(z)\frac{z\cdot\nabla\varphi(z)}{|z|^2}\;dz\\\nonumber\\
&=&(1+\frac{\alpha}{k})^2\left(\frac{a}{\pi k}-\frac{b}{\pi
A\sqrt{k}}\right),
\label{c}
\end{eqnarray}
where the constants $a$ and $b$ are given by
$$
a=2\dint_{1\leq|z|\leq
2}\;\log|z|\varphi(z)\frac{z\cdot\nabla\varphi(z)}{|z|^2}\;dz\quad\mbox{and}\quad
b=2\dint_{1\leq|z|\leq
2}\;\varphi(z)\frac{z\cdot\nabla\varphi(z)}{|z|^2}\;dz
$$
Finally, \eqref{I}, \eqref{II} together with \eqref{a}, \eqref{b} and \eqref{c} imply that for every $A>0$,
$$
 1+\frac{C_1}{\pi
A^2}\;\leq\;\liminf_{k\rightarrow\infty}\;H(g).
$$
\end{proof}
The main result of this section reads.
\begin{thm}
\label{tIPNLS} Let $\alpha>0$ and $A>0$ be real numbers, and
\begin{eqnarray*}
u_k(0,x)&=&g_{\alpha, A, k}(\nu_k(A)x),\\\\
v_k(0,x)&=&g_{0, A, k}(\nu_k(A)x).
\end{eqnarray*}
Denote by $u_k$ and $v_k$ the associated solutions of \eqref{eq1}. Then, there exists a
sequence $t_k\longrightarrow 0^+$ such that
\begin{equation}
\label{eqIPNLS}
\dliminf_{k\rightarrow\infty}\|\nabla(u_k-v_k)(t_k,\cdot)\|_{L^2(\R^2)}\gtrsim
1.
\end{equation}
\end{thm}

A general strategy to prove such instability result is to analyze
the associated ordinary differential equation
(see for instance, \cite{CCT, CCT1}). More precisely, let $\Phi$ solve
\begin{equation}
\label{ODE.NLS} \left\{
\begin{matrix}
i\;\partial_t \Phi(t,y)&=& f(\Phi(t,y)), \\\\
\Phi(0,y) &=& g_{\alpha, A, k}(y).\\
\end{matrix}
\right.
\end{equation}
The problem \eqref{ODE.NLS} has an explicit solution given by:
\begin{eqnarray*}
\label{ODE.SOLUTION}
 \Phi_0^{(\alpha,A,k)}(t,y)&=&g_{\alpha, A,
k}(y)\;\exp\left(-it(e^{4\pi g_{\alpha, A, k}(y)^2}
-1)\right)\\&:=&g_{\alpha, A, k}(y)\;\exp\left(-itK(g_{\alpha, A,
k})(y)\right)
\end{eqnarray*}
where $K(z)=e^{4\pi |z|^2} -1$.\\

In the case of a power type nonlinearity, the common element in all
arguments is a quantitative analysis of the NLS equation in the {\it small dispersion} limit
$$
i\partial_t \Phi+\nu^2\;\Delta \Phi=\sigma\;|\Phi|^{p-1}\Phi
$$
where (the dispersion coefficient) $\nu$ is small. Formally, as $\nu\rightarrow 0$ this
equation approaches the ODE
$$
i\partial_t \Phi=\sigma\;|\Phi|^{p-1}\Phi
$$
which has an explicit solution (see \cite{CCT, CCT1} for more
details). This fact and the invariance of equations of the type (\ref{NLSp}) under the
scaling $\Phi\mapsto\Phi^{\lambda}$ defined by
$$
\Phi^{\lambda}(t,x):=\lambda^{2/(1-p)}\;\Phi(\lambda^{-2}
t,\lambda^{-1}x)$$ play a crucial role in the ill-posedness results
obtained in \cite{CCT, CCT1} to make the decoherence happen during the approximation.

Unfortunately, no scaling leaves our
equation invariant and this seems to be the major difficulty since
it forces us to suitably construct the initial data in
Theorem \ref{tIPNLS}. Our solution to this difficulty (and others) proceeds in the context
of energy and Strichartz estimates for the following equation
\begin{equation}
\label{SC-NLS} i\partial_t \Phi+\nu^2\;\Delta \Phi=f(\Phi)
\end{equation}
It turns out that given the scale $ \nu_k(A)$, then for times close to $\frac{e^{-k}}{\sqrt{k}}$,
equation \eqref{SC-NLS} approaches the associated ODE
\eqref{ODE.NLS}.

\begin{proof}[Proof of Theorem \ref{tIPNLS}]
The proof is divided into two steps.\\

{\bf First Step:} \underline{\it"Decoherence"} \\

The key Lemma is the following.
\begin{lemma}
\label{L.IPNLS} Let ${\mathcal C}_k$ denote the ring
$\{2e^{-k/2}\leq|y|\leq 3e^{-k/2}\}$ and
$t_k^\varepsilon=\varepsilon\frac{e^{-k}}{\sqrt{k}}$. Then, we have
\begin{equation}
\label{eq.IPNLS} \varepsilon
e^{2\alpha}e^{-C\frac{\alpha}{k}}e^{-C\frac{\alpha^2}{k^2}}\lesssim
\|\nabla\Phi_0^{\alpha,A,k}(t_k^\varepsilon)\|_{L^2({\mathcal
C}_k)}\lesssim(1+\frac{\alpha}{k})^3\left(\frac{1}{k}+\varepsilon
e^{2\alpha} e^{C\frac{\alpha^2}{k}}\right)
\end{equation}
where $C$ stands for an absolute positive constant which may change
from term to term.
\end{lemma}
\begin{proof}[Proof of Lemma \ref{L.IPNLS}]
Write $\Phi_0^\alpha$ for $\Phi_0^{\alpha,A,k}$, then
$$
\|\nabla\Phi_0^\alpha(t)\|_{L^2({\mathcal C}_k)}^2=\|\nabla
g\|_{L^2({\mathcal C}_k)}^2+64\pi^2 t^2\|g^2 e^{4\pi g^2}\nabla
g\|_{L^2({\mathcal C}_k)}^2
$$

In view of the definition of $\eta_k$ and $\varphi$, we get
$$
\|\nabla\Phi_0^\alpha(t)\|_{L^2({\mathcal
C}_k)}^2=(1+\frac{\alpha}{k})^2 {\mathcal I}+64\pi^2
t^2(1+\frac{\alpha}{k})^6 {\mathcal J}
$$
where \begin{eqnarray*} {\mathcal I}&=&\dint_{{\mathcal C}_k}|\nabla
f_k(y)|^2dy=\frac{2}{k}\log(\frac{3}{2})\;\lesssim \frac{1}{k}\\\\
{\mathcal J}&=&\dint_{{\mathcal C}_k}\;|f_k(y)|^4|\nabla
f_k(y)|^2e^{8\pi(1+\frac{\alpha}{k})^2f_k(y)^2}\;dy\\\\
&=&\frac{2}{\pi^2 k^2}\dint_{2e^{-k/2}}^{3e^{-k/2}}\;
\exp\left({\frac{8}{k}(1+\frac{\alpha}{k})^2\log^2 r}\right)\log^4
r\;\frac{dr}{r}
\end{eqnarray*}
We conclude the proof by remarking that, for $2 e^{-k/2}\leq r\leq 3
e^{-k/2}$, we have
$$
k^4 e^{k/2} e^{4\alpha}
e^{2k}e^{-C\frac{\alpha}{k}}e^{-C\frac{\alpha^2}{k^2}}\lesssim
\exp\left({\frac{8}{k}(1+\frac{\alpha}{k})^2\log^2
r}\right)\frac{\log^4 r}{r}\lesssim k^4 e^{k/2} e^{4\alpha}
e^{2k}e^{C\frac{\alpha}{k}}
$$
\end{proof}

\begin{cor}
Let $\alpha> 0$ be a real number. Then,
\begin{equation} \label{decho}
\dliminf_{k\rightarrow\infty}\Big\|\nabla \left(\Phi_0^{\alpha,
A,k}-\Phi_0^{0,A,
k}\right)(t_k^\varepsilon)\Big\|_{L^2(\R^2)}\gtrsim\varepsilon(e^{2\alpha}-1).
\end{equation}
\end{cor}
\begin{proof}
In view of the previous lemma, we have
\begin{eqnarray*}
\Big\|\nabla\left(\Phi_0^\alpha-\Phi_0^0\right)(t_k^\varepsilon)\Big\|_{L^2(\R^2)}&\geq&
\Big\|\nabla\left(\Phi_0^\alpha-\Phi_0^0\right)(t_k^\varepsilon)\Big\|_{L^2({\mathcal
C}_k)}\\\\&\geq&\Big\|\nabla\Phi_0^\alpha(t_k^\varepsilon)\Big\|_{L^2({\mathcal
C}_k )}- \Big\|\nabla\Phi_0^0(t_k^\varepsilon)\Big\|_{L^2({\mathcal
C}_k )}\\\\&\gtrsim&\varepsilon
e^{2\alpha}e^{-C\frac{\alpha}{k}}e^{-C\frac{\alpha^2}{k^2}}-\left(\frac{1}{k}+\varepsilon\right)
\end{eqnarray*}
and the conclusion follows.
\end{proof}

{\bf Second Step:} \underline{\it Approximation} \\
The end of the proof of Theorem \ref{tIPNLS} lies in the
following technical lemmas.
\begin{lemma}
\label{L.ODEH1} The solution $\Phi^{\alpha, A, k}_0$ of
\eqref{ODE.NLS} satisfies
\begin{equation}
\label{estH2.ODE}
\|\nabla^3\Phi^{\alpha,A,k}_0(t)\|_{L^2}\lesssim\;\frac{e^k}{\sqrt{k}}\left(1+tk^{1/3}\;
e^k\right)^3
\end{equation}
\end{lemma}
\begin{proof}[Proof of Lemma \ref{L.ODEH1}]\quad\\
Write $\Phi_0$ for $\Phi^{\alpha,A,k}_0$ and $g$ for
$g_{\alpha,A,k}$. Clearly,
\begin{eqnarray*}
\nabla\Phi_0&=&\Big(\nabla g-it g K'(g)\nabla
g\Big)e^{-itK(g)}:=g_1e^{-itK(g)},\\\\
\nabla^2\Phi_0&=&\Big(\nabla
g_1-itg_1K'(g)\nabla g\Big)e^{-itK(g)}:=g_2e^{-itK(g)},\\\\
\nabla^3\Phi_0&=&\Big(\nabla g_2-itg_2K'(g)\nabla
g\Big)e^{-itK(g)}:=g_3e^{-itK(g)},
\end{eqnarray*}
so
\begin{eqnarray*}
\|\nabla^3\Phi_0\|_{L^2}&\lesssim&\|\nabla
g_2\|_{L^2}+t\|g_2 K'(g)\nabla g\|_{L^2},\\\\
&\lesssim&\|\nabla^3 g\|_{L^2}+tA_1+t^2A_2 + t^3 A_3,
\end{eqnarray*}
where
\begin{eqnarray*}
A_1&=&\|K'(g)\nabla^2g\nabla g\|_{L^2}+\|K^{''}(g)(\nabla
g)^3\|_{L^2}\\\\&+&\|g K^{'''}(g)(\nabla g)^3\|_{L^2}+
\|g K^{''}(g)\nabla^2 g\nabla g\|_{L^2}+\|g K'(g)\nabla^3 g\|_{L^2},\\\\
A_2&=&\|(K'(g))^2(\nabla g)^3\|_{L^2}+\|g K'(g) K^{''}(g)(\nabla
g)^3\|_{L^2}+\|g (K'(g))^2\nabla^2 g\nabla g\|_{L^2},\\\\
A_3&=&\|g(K'(g))^3(\nabla g)^3\|_{L^2}.
\end{eqnarray*}
Now,
\begin{eqnarray*}
\|\nabla^3
g\|_{L^2}^2&\lesssim&\dint_{\frac{3}{2}e^{-k/2}\leq|y|\leq
2\nu_k(A)} \;|\nabla^3 f_k|^2 dy+
\mbox{l.o.t},\\\\
&\lesssim&\frac{1}{k}\dint_{\frac{3}{2}e^{-k/2}}^{2\nu_k(A)}\;\frac{dr}{r^5}+\mbox{l.o.t},\\\\
&\lesssim&\frac{e^{2k}}{k}.
\end{eqnarray*}
On the other hand
\begin{eqnarray*}
\|g(K'(g))^3(\nabla g)^3\|_{L^2}^2&\lesssim&\dint\;|g|^8 e^{24\pi
g^2}\;|\nabla g|^6\;dy\\\\
&\lesssim&\frac{1}{k^7}\dint_{\frac{3}{2}e^{-k/2}}^{2\nu_k(A)}\;\log^8
r \;e^{\frac{24}{k}(1+\frac{\alpha}{k})^2\log^2
r}\;\frac{dr}{r^5}+\mbox{l.o.t}
\end{eqnarray*}
\end{proof}
The next lemma states the energy and
Strichartz estimates for NLS with small dispersion coefficient
$\nu$.
\begin{lemma}
\label{L.ENSTRI}
 Let $v_0$ be a function in $H^1(\R^2)$ and $F\in
L^1(\R, H^1(\R^2))$. Denote by $v$ the solution of the inhomogeneous
linear Schr\"odinger equation
$$
i\partial_t v+\nu^2\;\Delta_y v=F(t,y)
$$
with initial data $u(0,y)=v_0(y)$. \\
Then, a constant $C$ exists such that for any $T>0$, we have

\begin{eqnarray}
\label{eq.ENSTRI} \nonumber \|\nabla
v\|_{L^\infty_T(L^2)}+\frac{1}{\nu}\|
v\|_{L^\infty_T(L^2)}+\nu^{1/2}\|v\|_{L^4_T(\dot{\mathcal
C}^{1/2})}&\lesssim&\|\nabla
v_0\|_{L^2}+\frac{1}{\nu}\|v_0\|_{L^2} \\
 &+&\|\nabla F\|_{L^1_T(L^2)}+\frac{1}{\nu}\|F\|_{L^1_T(L^2)}.
\end{eqnarray}
\end{lemma}

This lemma can be obtain from the standard Strichartz estimates through an obvious scaling. It can be seen as a semiclassical Strichartz estimate
which permits an extension of the approximation time. Also, this
lemma plays a role in the NLS analysis that is played by finite
propagation speed in the corresponding NLW arguments.
Now we are ready to end the proof of Theorem \ref{tIPNLS}. For this purpose, denote (for
simplicity) by $w:= \Phi-\Phi_0$ where $\Phi_0$ is given by
\eqref{ODE.SOLUTION} and $\Phi$ solves the problem

\begin{equation}
\label{eqNU} \left\{
\begin{matrix}
i\partial_t \Phi(t,y)+\nu^2\;\Delta_y\Phi(t,y)= f(\Phi(t,y)),\\
\Phi(0,y) = g(y).\\
\end{matrix}
\right.
\end{equation}
Set
\begin{eqnarray*}
M_0(w,t)&{\buildrel\hbox{\footnotesize def}\over =}&\|
w\|_{L^\infty(([0,t];L^2)}\\
M_1(w,t)&{\buildrel\hbox{\footnotesize def}\over =}&\|\nabla
w\|_{L^\infty(([0,t];L^2)}+\frac{1}{\nu}\|
w\|_{L^\infty(([0,t];L^2)}+\nu^{1/2}\|w\|_{L^4(([0,t];\dot{\mathcal
C}^{1/2})}
\end{eqnarray*}
We will prove the following result.
\begin{lemma}
\label{Approximation} For
$t_k^\varepsilon\thickapprox\varepsilon\frac{e^{-k}}{\sqrt{k}}$ and
$k$ large, we have
\begin{eqnarray}
\label{Approx1} M_0(w,t_k^\varepsilon)\lesssim\;e^{-k/2}\;\nu^{3/2}.
\end{eqnarray}
\begin{eqnarray}
\label{Approx2} M_1(w,t_k^\varepsilon)\lesssim\;\nu^{1/2}.
\end{eqnarray}
\end{lemma}
\begin{proof}[Proof of Lemma \ref{Approximation}]
Since $w$ solves
$$
i\partial_t w+\nu^2\;\Delta_y w=
f(\Phi_0+w)-f(\Phi_0)-\nu^2\;\Delta_y\Phi_0,\quad w(0,y)=0,
$$
then using the $L^2$ energy estimate, we have
$$
M_0(w,t)\lesssim\;\nu\;I_2(t)+\nu\;I_4(t),
$$
where we set
\begin{eqnarray*}
I_2(t):&=&\frac{1}{\nu}\;\|f(\Phi_0+w)-f(\Phi_0)\|_{L^1(([0,t];L^2)},\\\\
I_4(t):&=&\nu\;\|\nabla^2\Phi_0\|_{L^1(([0,t];L^2)}.
\end{eqnarray*}
Note that we have the following
\begin{eqnarray*}
\|\big(f(\Phi_0+w)-f(\Phi_0)\big)(t,\cdot)\|_{L^2}&\lesssim&\|w(t,\cdot)\|_{L^2}
\Big(\|\big(\Phi_0+w\big)(t,\cdot)\|_{L^\infty}^2\;e^{4\pi\|\big(\Phi_0+w\big)(t,\cdot)\|_{L^\infty}^2}\\&+&
\|w(t,\cdot)\|_{L^\infty}^2\;e^{4\pi\|w(t,\cdot)\|_{L^\infty}^2}\Big).
\end{eqnarray*}
Integrating in time we have
\begin{eqnarray*}
\|\big(f(\Phi_0+w)-f(\Phi_0)\big)(t,\cdot)\|_{L^1([0,T],L^2)}&\lesssim& \int_0^tM_0(w,s)\Big(\|\big(\Phi_0+w\big)(s,\cdot)\|_{L^\infty}^2\;e^{4\pi\|\big(\Phi_0+w\big)(s,\cdot)\|_{L^\infty}^2}\\&+&
\|w(s,\cdot)\|_{L^\infty}^2\;e^{4\pi\|w(s,\cdot)\|_{L^\infty}^2}\Big)ds.
\end{eqnarray*}
Using Lemma \ref{L.ODEH1} and the following simple fact
\begin{eqnarray}
\label{exercice}
\sup_{x\geq 0}\left(x^m\;e^{-\gamma
x^2}\right)=\left(\frac{m}{2\gamma}\right)^{m/2}\;e^{-m/2},\quad m\in\N,\quad \gamma>0,
\end{eqnarray}
we deduce that
\begin{equation}
\label{Gronw-0} M_0(w,t)\lesssim h_0(t)+\dint_0^t\;A_0(s)M_0(w,s)ds,
\end{equation}
where we set

\begin{eqnarray*}
h_0(t)&=&\nu^2\;e^{k/2}\dint_0^t\;(1+ske^k)^2\;ds\lesssim t\nu^2\;e^{k/2}(1+(tke^k)^2),\\\\
A_0(s)&=&k\;e^{4\pi(1+1/k)\|\big(\Phi_0+w\big)(s,\cdot)\|_\infty^2}+k\;
e^{4\pi(1+1/k)\|w(s,\cdot)\|_\infty^2}.
\end{eqnarray*}
Applying the logarithmic inequality \eqref{eq6} ( for $\lambda=\frac{1}{\pi}$) and using the fact that
$M_1(w,t)\lesssim\nu^{1/2}$, we obtain

$$
A_0(s)\lesssim
ke^k\left(C+\nu^{-1/2}\;\|w\|_{C^{1/2}}\right)^{4\delta_1(k)}+k\;
\left(C+\nu^{-1/2}\;\|w\|_{C^{1/2}}\right)^{4\delta_2(k)}
$$
where

$$
\delta_1(k)=(1+1/k)(\nu+\nu^{1/2}\;\sqrt{k})\quad\mbox{and}\quad
\delta_2(k)=(1+1/k)\nu.
$$
Now, using H\"older inequality in time we deduce that

\begin{eqnarray*}
\dint_0^t\;A_0(s)ds&\lesssim&ke^k\;t^{1-\delta_1(k)}\left(t^{1/4}+
\nu^{-1/2}\;\|w\|_{L^4_t(C^{1/2})}\right)^{4\delta_1(k)}\\\\&+&k\;t^{1-\delta_2(k)}\left(t^{1/4}+
\nu^{-1/2}\;\|w\|_{L^4_t(C^{1/2})}\right)^{4\delta_2(k)}\\\\
&\lesssim&k e^k\;t^{1-\delta_1(k)}\left(t^{1/4}+
\frac{1}{\sqrt{\nu}}\right)^{4\delta_1(k)}+\nu\;k^{3/2}\;
t^{1-\delta_2(k)}\left(t^{1/4}+
\frac{1}{\sqrt{\nu}}\right)^{4\delta_2(k)}
\end{eqnarray*}
It is easy to see that for $t\approx t_k^\varepsilon$,
$$
\dint_0^t\;A_0(s)ds\lesssim \varepsilon\sqrt{k}.
$$
Hence, by Gronwall's lemma

\begin{eqnarray*}
M_0(w,t_k^\varepsilon)&\lesssim& h_0(t_k^\varepsilon)\exp(C\varepsilon\sqrt{k})\\\\
&\lesssim&\nu^{3/2}\;e^{-k/2}\;\frac{\varepsilon}{\sqrt{k}}\;(1+\varepsilon^2\;k)\exp\left(C\varepsilon\sqrt{k}-\frac{\sqrt{k}}{2A}\right)
\\\\
&\lesssim&\nu^{3/2}\;e^{-k/2}
\end{eqnarray*}
provided that $\varepsilon<\frac{c}{A}$.

Similarly we proceed for $M_1$. According to Lemma \ref{L.ENSTRI}, we have
$$
 M_1(w,t)\lesssim I_1(t)+I_2(t)+I_3(t)+I_4(t),
$$
where in addition we set
\begin{eqnarray*}
I_1(t):&=&\|\nabla(f(\Phi_0+w)-f(\Phi_0))\|_{L^1(([0,t];L^2)},
\end{eqnarray*}
and
\begin{eqnarray*}
I_3(t):&=&\nu^2\;\|\nabla^3\Phi_0\|_{L^1(([0,t];L^2)}.
\end{eqnarray*}
Note that

\begin{eqnarray*}
\noindent\|\nabla\left(f(\Phi_0+w)-f(\Phi_0)\right)\|_{L^2_x}&\lesssim&\|w\|_{L^2_x}
\|\nabla\Phi_0\|_{L^\infty_x}\Big(\|\Phi_0\|_{L^\infty_x}\;e^{4\pi\|\Phi_0\|_{L^\infty_x}^2}\\&+&
\|\Phi_0\|_{L^\infty_x}^2\;e^{4\pi\|\Phi_0\|_{L^\infty_x}^2}+
\|w\|_{L^\infty_x}\;e^{4\pi\|w\|_{L^\infty_x}^2}\\&+&
\|w\|_{L^\infty_x}^3\;e^{4\pi\|w\|_{L^\infty_x}^2}\Big)\\&+& \|\nabla
w\|_{L^2_x}\|\Phi_0+w\|_{L^\infty_x}^2\;e^{4\pi\|\Phi_0+w\|_{L^\infty_x}^2}\\&+&
\|\nabla w\|_{L^2_x}e^{4\pi\|\Phi_0+w\|_{L^\infty_x}^2}.
\end{eqnarray*}
Arguing as before, we have

\begin{eqnarray}
\label{Gronw} M_1(w,t)&\lesssim& h_1(t)+\frac{1}{\nu}\dint_0^t\;A_0(s)M_0(w,s)ds+\dint_0^t\;A_1(s)M_0(w,s)ds\\ \nonumber&+&
\dint_0^t\;A_0(s)M_1(w,s)ds
\end{eqnarray}
with in addition
\begin{eqnarray*}
h_1(t)&=&\nu\;\|\nabla^3\Phi_0\|_{L^1_t(L^2_x)}\lesssim\;\frac{e^k}{\sqrt{k}}\left(1+tk^{1/3}\;
e^k\right)^3,\\\\
A_1(s)&=&\sqrt{k}\;e^{k/2}(1+s\;k\;e^k)(e^k+\sqrt{k}\;e^{4\pi(1+1/k)\|w\|_\infty^2}).
\end{eqnarray*}

Here we used the Poincar\'e inequality and \eqref{exercice}. Now we return to $M_1(w,t)$ for which we have to prove (\ref{Approx2}). Using Lemma \ref{L.ODEH1}, (\ref{Approx1}) and Logarithmic
inequality, we get for
$t\approx t_k^\varepsilon$
\begin{eqnarray*}
 h_1(t)&\lesssim&\nu\;\frac{\varepsilon}{k}\;\left(1+\frac{\varepsilon^3}{\sqrt{k}}\right),\\\\
\frac{1}{\nu}\;\dint_0^t\;M_0(w,s)\;A_0(s)\;ds&\lesssim&\;e^{-k/2}\;\nu^{1/2}\;\varepsilon\;\sqrt{k},\\\\
\dint_0^t\;M_0(w,s)\;A_1(s)\;ds&\lesssim&\;\nu^{3/2}\;\varepsilon\;k\;(1+\varepsilon\;\sqrt{k}).
\end{eqnarray*}

Gronwall's lemma yields
$$
M_1(w,t)\lesssim\;\nu^{1/2}\Big(\nu^{1/2}+\varepsilon\;\sqrt{k}\;e^{-k/2}+\nu\;\varepsilon\;k\;(1+\varepsilon\sqrt{k})\Big)\exp(C\varepsilon\sqrt{k})\lesssim\;\nu^{1/2}
$$
provided that $\varepsilon<\frac{c}{A}$. This completes the proof of
Lemma \ref{Approximation}
\end{proof}

Finally, a comparison of \eqref{eq.IPNLS} with the approximation bounds
\eqref{Approx1}, \eqref{Approx2} implies \eqref{eqIPNLS}. This
completes the proof of Theorem \ref{tIPNLS}.

\end{proof}
%@@@@@@@@@@@@@@@@@@@@@@@@@@@@@@@@@@@@%@@@@@@@@@@@@@@@@@@@@@@@@@@@@@@@@@@@@%@

\end{document}